%% file: 2003-8.tex
\documentclass{gtart}
\input gtoutput

\lognumber{236}
\volumenumber{7}\papernumber{8}\volumeyear{2003}
\pagenumbers{287}{309}
\received{13 February 2002}
\revised{18 March 2003}
\accepted{24 April 2003}
\published{27 April 2003}
\proposed{David Gabai}
\seconded{Cameron Gordon, Joan Birman}

\usepackage{amsmath,amssymb}
\input rlepsf

\newtheorem{theorem}{Theorem}[section]
\newtheorem{lemma}[theorem]{Lemma}

\newtheorem{proposition}[theorem]{Proposition}

\theoremstyle{remark}
\newtheorem{remark}[theorem]{Remark}
\newtheorem{definition}[theorem]{Definition}
\newtheorem*{acknowledgments}{Acknowledgments}

\newcommand{\ZZ}{\mathbb{Z}}

\newcommand{\RR}{\mathbb{R}}

\def\g{\gamma}

\begin{document}

\title{An algorithm to detect laminar 3-manifolds}
\author{Ian Agol\\Tao Li}
\asciiaddress{Department of Mathematics, University of Illinois at
Chicago\\322 SEO m/c 249, 851 S. Morgan Street\\Chicago, IL
60607-7045, USA\\Department of Mathematics, 
401 Math Sciences\\Oklahoma State University\\Stillwater, OK 74078-1058, USA}

\address{Department of Mathematics, University of Illinois at
Chicago\\322 SEO m/c 249, 851 S. Morgan Street\\Chicago, IL
60607-7045, USA}
\secondaddress{Department of Mathematics, 
401 Math Sciences\\Oklahoma State University\\Stillwater, OK 74078-1058, USA}

\asciiemail{agol@math.uic.edu, tli@math.okstate.edu}

\email{agol@math.uic.edu}
\secondemail{tli@math.okstate.edu}

\url{http://www.math.uic.edu/\char'176agol,
http://www.math.okstate.edu/\char'176tli/}

\begin{abstract}
We show that there are algorithms to determine if a 3-manifold
contains an essential lamination or a Reebless foliation.
\end{abstract}


\maketitlepage

\section{Introduction} \label{intro}

Essential laminations were introduced by Gabai and Oertel
\cite{GO}, as a generalization of incompressible surfaces
\cite{Hak68}, measured incompressible laminations \cite{MS88},
Reebless foliations \cite{No65}, and laminations coming from
pseudo-Anosov flows of fibred manifolds \cite{Th88}. Many
properties of Haken manifolds are now known to hold for manifolds
containing essential laminations, eg, they satisfy a weak
hyperbolization property: either the fundamental group contains
$\ZZ^2$, or the manifold has word-hyperbolic fundamental group, as
shown by Gabai and Kazez for genuine essential laminations
\cite{GK}, and by Calegari for taut foliations \cite{C,C3,C02}. By
the solution of the Seifert conjecture \cite{CJ,Ga9}, if
$\ZZ^2\leq \pi_1(M)$, then either $M$ contains an embedded
incompressible torus, or it is a small Seifert fibered space with
infinite fundamental group, for which essential laminations have
been classified \cite{Br2,Br3,EHN,JN,Na}. Thus, it is natural to
wonder to what extent 3-manifolds with infinite fundamental group
and no $\ZZ^2$ subgroup of the fundamental group contain essential
laminations.

To resolve this question, it would be useful to have an algorithm
to determine whether or not a 3-manifold contains an essential
lamination. Jaco and Oertel proved that there is an algorithm to
determine if an irreducible 3-manifold is Haken, by showing that
if there is an incompressible surface, then there is one among a
finite collection of algorithmically constructible normal
surfaces. Oertel showed that there is an algorithm to determine
whether a manifold contains an affine lamination \cite{O3}. Some
progress was made by Brittenham, who showed that if a manifold
contains an essential lamination, then it contains an essential
lamination which is normal with respect to a given triangulation
\cite{Br}. Brittenham's normalization process was analyzed in
depth by Gabai \cite{Ga8}, who determined exactly how the
normalized lamination differs from the starting lamination. Li
gave a criterion on a branched surface embedded in a 3-manifold,
called a {\it laminar branched surface}, which implies that the
manifold has an essential lamination carried by the branched
surface, and he proved that any manifold with an essential
lamination contains such a laminar branched surface \cite{Li}. We
show in this paper that the results of Gabai and Li imply the
existence of algorithms to determine if a 3-manifold contains an
essential lamination or a Reebless foliation, answering Problem
2.3 of \cite{Ga7}. Recently, Roberts, Shareshian, and Stein
\cite{RSS} have given examples of 3-manifolds which admit no
Reebless foliation, and Sergio Fenley \cite{F02} has extended this
to show that their examples contain no essential lamination. The
method of proof they use is to show that the groups do not act on
order trees. It is unknown whether there are non-laminar
3-manifolds that do not act on order trees, so their approach does
not provide an algorithm.

The algorithm to determine if a manifold has an essential
lamination described in theorem \ref{alg} proceeds by first
finding a finite collection of incompressible Reebless branched
surfaces which carry every nowhere dense essential lamination in
$M$ which is normal with respect to a particular one-efficient
triangulation. Then the algorithm splits into two procedures run
in tandem: for each branched surface, one procedure tries to find
a laminar branched surface carried by it, which succeeds for at
least one branched surface if the manifold contains an essential
lamination; the other procedure tries to split each branched
surface in all possible ways, which fails for every branched
surface in the finite collection if the manifold does not contain
an essential lamination.

\begin{acknowledgments}
We would like to thank Bus Jaco for many helpful conversations
about one-efficient triangulations, and Mark Brittenham for
finding some errors in an earlier draft of this paper. We also
thank the referee for many useful comments, and for noticing some
gaps and errors in earlier drafts. The second author is partially
supported by NSF grant DMS 0102316.
\end{acknowledgments}

\section{Laminar branched surfaces}

In this paper, we will assume that all 3-manifolds are orientable,
since irreducible non-orientable manifolds are Haken, and
therefore laminar.

\begin{figure}
\cl{\relabelbox\small
\epsfxsize=4.0in\epsfbox{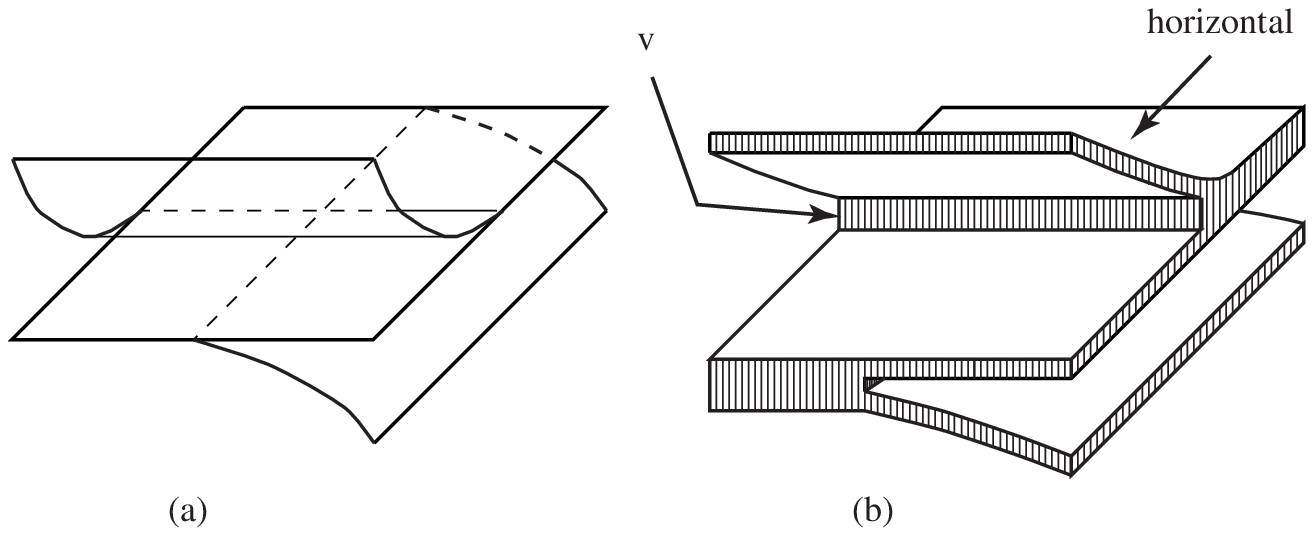}
\relabel {(a)}{(a)} 
\relabel {(b)}{(b)}
\relabel {horizontal}{$\partial_hN(B)$} 
\relabel {v}{$\partial_vN(B)$}
\endrelabelbox}
\nocolon\caption{}\label{branched}
\end{figure}

\begin{definition}
A {\it branched surface}  $B$ is a union of finitely many compact
smooth surfaces glued together to form a compact subspace (of $M$)
locally modelled on Figure \ref{branched}(a).
\end{definition}

{\bf Notation}\qua Throughout this paper, we denote the interior of
$X$ by $int(X)$, and denote the number of components of $X$ by
$|X|$, for any $X$.\medskip

Given a branched surface $B$ embedded in a 3-manifold $M$, we
denote by $N(B)$ a regular neighborhood of $B$, as shown in Figure
\ref{branched}(b). One can regard $N(B)$ as an interval bundle
over $B$, with a foliation $\mathcal{V}$ by intervals. We denote
by $\pi:N(B)\to B$ the projection that collapses every interval
fiber of $\mathcal{V}$ to a point. The {\it branched locus} of $B$
is the 1-skeleton of $B$, when viewed as a spine. So, $L$ can be
considered as a union of immersed curves in $B$, and we call a
point in $L$ a {\it double point} of $L$ if any small neighborhood
of this point is modelled on Figure \ref{branched}(a).

A (codimension one) lamination $\lambda$ in a 3-manifold $M$ is a
foliated, closed subset of $M$, ie, $\lambda$ is covered by
a collection of open sets $U$ of the form $\RR^2\times \RR$ such
that $\lambda\cap U = \RR^2\times C$, where $C$ is a closed set in
$\RR$, and the transition maps preserve the product structures.
The coordinate neighborhoods of leaves are of the form
$\RR^2\times x $, $x\in C$. Let $M_\lambda$ be the metric
completion of the manifold $M-\lambda$ with the path metric
inherited from a metric on $M$. Let $H=\{(x,y)\in\RR^2|y\geq 0\}$
be the closed upper half plane. An {\it end compression} is a
proper embedding $d:(H,\partial H)\to (M_\lambda,\partial
M_\lambda)$ such that $d|_{\partial H}$ does not extend to a
proper embedding $d':H\to
\partial M_\lambda$.

Given a branched surface $B\subset M$, a lamination $\lambda$ is
{\it carried} by $B$ if $\lambda\subset N(B)$ and each leaf of
$\lambda$ is transverse to $\mathcal{V}$, and {\it fully carries}
if $\lambda$ meets every fiber of $\mathcal{V}$. Similarly, a
branched surface $B'$ is carried by $B$ if $B'\subset N(B)$ and
$B'$ is smoothly transverse to $\mathcal{V}$. If $B'$ is carried
by $B$ and meets every fiber of $\mathcal{V}$, then $B'$ is a {\it
splitting} of $B$.

\begin{definition}
$\lambda$ is an {\it essential lamination} in $M$ if it satisfies the
following conditions:
\begin{enumerate}
\item
The inclusion of leaves of the lamination into $M$ induces a injection
on $\pi_1$.
\item
$M_\lambda$ is irreducible.
\item
$\lambda$ has no sphere leaves.
\item
$\lambda$ has no end compressions.
\end{enumerate}
\end{definition}

\begin{theorem}{\rm\cite[Proposition 4.5]{GO}} \label{lamdef}

\begin{enumerate}

\item

Every essential lamination is fully carried by a branched surface
with the following properties:
\begin{enumerate}
\item[\rm(a)]
$\partial_h N(B)$ is incompressible in $M-int(N(B))$, no component
of $\partial_h N(B)$ is a sphere, and $M-B$ is irreducible.
\item[\rm(b)]
There is no monogon in $M-int(N(B))$, ie, no disk $D\subset
M-int(N(B))$ with $\partial D =D\cap N(B)=\alpha\cup\beta$, where
$\alpha\subset\partial_v N(B)$ is in an interval fiber of
$\partial_v N(B)$ and $\beta\subset\partial_h N(B)$.
\item[\rm(c)]
There is no Reeb component, ie, $B$ does not carry a
sublamination of a Reeb foliation of a solid torus.
\item[\rm(d)]
$B$ has no disk of contact, ie, no disk $D\subset N(B)$
such that $D$ is transverse to the $I$-fibers $\mathcal{V}$ of $N(B)$ and
$\partial D\subset \partial_v N(B)$, see Figure \ref{disk}(a) for an
example.

\begin{figure}
\cl{\relabelbox\small
\epsfxsize=4.0in\epsfbox{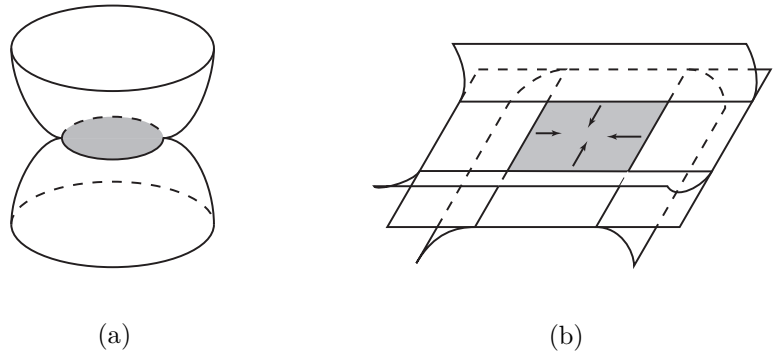}
\relabel {(a)}{(a)} 
\relabel {(b)}{(b)}
\endrelabelbox}
\caption{Sink disks}\label{disk}
\end{figure}

\end{enumerate}
\item
If a branched surface with the properties above fully carries a lamination,
then it is an essential lamination.
\end{enumerate}
\end{theorem}

A branched surface satisfying conditions (a), (b), and (d) of
Theorem \ref{lamdef} is an {\it incompressible branched surface},
and if it satisfies conditions (a)--(d), it is an {\it
incompressible Reebless branched surface}.

\begin{definition}
Suppose $B$ is an incompressible branched surface. Let $D_1$ and
$D_2$ be two disk components of $\partial_h N(B)$, such that
$\pi(\partial D_1)=\pi(\partial D_2)$. Then there is a region
$K\cong D^2\times I$ of $M-int(N(B))$ such that $D_1\cup D_2 =
D^2\times \partial I$.
  If $\pi(int (D_1))\cap\pi(int (D_2))=\emptyset$, then we call $K$ a
  {\it trivial bubble}, and we may collapse $B$ along $K$ by collapsing
  each $I$-fiber to a point, to get a new branched surface $B'$, such
  that any lamination is carried by $B$ if and only if it is carried by $B'$.
\end{definition}

\begin{definition}
A {\it sink disk} of $B$ is a disk component $D$ of $B-L$ such
that $\partial_h\pi^{-1}(\overline{D}) \cap
\partial_v N(B) =\emptyset$ (see figure \ref{disk}(b)).
\end{definition}

\begin{definition} \label{lambrdef}
A branched surface $B$ having no trivial bubble in $M$ is called a
{\it laminar branched surface} if it satisfies conditions 1(a)--(c)
of theorem \ref{lamdef} and has no sink disk.
\end{definition}

The following theorem follows from the proof of Lemma 5.4 in
\cite{Li}.

\begin{theorem} \label{li}
A branched surface carrying an essential lamination which is not a
lamination by planes has a splitting which is laminar.
\end{theorem}

\begin{theorem} \label{lamalg1}
Given a branched surface $B$ which carries no spheres or tori,
there is a procedure which terminates if and only if $B$ fully
carries a lamination. Moreover, if the procedure terminates, it
outputs a certificate that the branched surface fully carries a
lamination.
\end{theorem}
\begin{proof}

If $B$ is an essential branched surface, then by theorem \ref{li},
there is a laminar branched surface $B'$ fully carried by $B$,
ie, $B'$ has no sink disks and no trivial bubbles (one need
not check the other conditions for laminar branched surface, since
$B'$ cannot carry a sphere or torus, and $B'$ will have no monogon
or compression of $\partial_h N(B')$ in $N(B)$ since $B'$ is a
splitting of $B$, by Lemma 2.5 of \cite{GO}). Choose a
triangulation $\tau$ of $N(B)$ such that each simplex is
transverse to the fibers $\mathcal V$ of $N(B)$, and the foliation
of the simplex  by fibers of $\mathcal V$ is linear in the affine
structure on the simplex. There is a triangulation of $N(B)$
$\eta$ such that $B'$ is embedded in $\eta^{(2)}$. The
triangulations $\eta$ and $\tau$ have a common subdivision $\mu$,
which carries $B$ in $\mu^{(2)}$, since $\eta$ does. If we choose
$\mu$ to be generic, then it will be transverse to $\mathcal{V}$
since $\mathcal{V}$ is linear on each simplex of $\tau$. The
procedure proceeds in steps. Let $m=|\tau - \tau^{(2)}|$ be the
number of 3-simplices in $\tau$. Then the $n$th step of the
procedure first enumerates all subdivisions $\mu$ of $\tau$
transverse to $\mathcal{V}$ which have $m+n$ 3-simplices, up to
isotopy transverse to $\mathcal{V}$. For a given $n$, there are
only finitely many such subdivisions, which may be found
combinatorially, because of the linearity of the foliation
$\mathcal V$ restricted to each tetrahedron of $\tau$. For each
such subdivision $\mu$, it then enumerates all branched surfaces
in $\mu^{(2)}$ transverse to $\mathcal V$, and checks to see
whether they have no sink disk or trivial bubble. The certificate
is a description of the subdivision $\mu$ and the branched surface
in $\mu^{(2)}$.
\end{proof}

\begin{lemma} \label{incReeb}
For a 3-manifold $M$ with triangulation $\nu$, and
 a branched surface $B\subset \nu^{(2)}$, there is an
 algorithm to determine if $B$ is incompressible and
 Reebless.
\end{lemma}
\begin{proof}
We need to check that $B$ satisfies all the conditions of an
incompressible Reebless branched surface. The complement of $B$
has a triangulation $\nu'$ by splitting $\nu$ along $B$. Property
1(a) of theorem \ref{lamdef} states that $\partial_h N(B)$ is
incompressible in $M-int(N(B))$, no component of $\partial_h N(B)$
is a sphere, and $M-B$ is irreducible. To check that $\partial_h
N(B)$ is incompressible, we use Haken's algorithm, as improved by
Jaco and Tollefson in Algorithm 9.6 of \cite{JT}, by showing that
there is no vertex normal compressing disk for $\partial_h N(B)$.
One computes the Euler characteristic and number of boundary
components of each component of $\partial_h N(B)$ in order to
check that no component is a sphere. If $M-B$ is reducible, then
there would be an essential 2-sphere in $M-B$ which is a vertex
normal surface with respect to $\nu'$ by Theorem 5.2 \cite{JT}. We
may use an algorithm of Rubinstein \cite{Ru2} and Thompson
\cite{Tho} to show that every vertex normal 2-sphere in $M-B$
bounds a ball.

Property 1(b) of theorem \ref{lamdef} states that there is no
monogon in $M-int(N(B))$. Again, one need only check that there is
no vertex normal disk meeting $\partial_h N(B)$ exactly once using
the argument in Theorem 6.1 of \cite{JT}.

Property 1(c) states that there is no Reeb component, ie,
$B$ does not carry a torus that bounds a solid torus in $M$. To
check this, we compute the normal solution space of surfaces
carried by $B$, and compute all vertex surfaces which are tori.
For each one, we check that it does not bound a solid torus, by
checking that either it bounds no compressing disk using Haken's
algorithm, or that it has a compressing disk and bounds a ball
with knotted hole, which may be checked using  Rubinstein and
Thompson's algorithm to recognize the 3-ball \cite{Ru2,Tho}.

Property 1(d) states that there is no disk of contact. We may
solve relative normal coordinates to check that no surface carried
by $B$ with boundary on $\partial_v N(B)$ is a disk, using the
method of Floyd and Oertel in Claim 3, page 123 of \cite{FO}.
\end{proof}

\section{Splitting branched surfaces}
Given a branched surface $B\subset M$, let the pared locus
$p(B)\subset \partial N_v(B)$ be a collection of simple closed
curves which are cores for $\partial_v N(B)$.

\begin{definition}
A {\it splitting surface} $S\subset N(B)$ is a complete $1-1$
immersed surface carried by $N(B)$ such that $\partial S = p(B)$,
and such that for each $I$-fiber $J$ of $\mathcal{V}$, $S\cap J$
is a discrete set. There is a regular open neighborhood
$\mathcal{N}(S)$ which fibers over $S$ induced by the fibration
$\mathcal{V}$, such that $\mathcal{N}(S)\cap \partial_v
N(B)=int(\partial_v N(B))$ and $\partial
\overline{\mathcal{N}(S)}$ is an immersed surface transverse to
$\mathcal{V}$.
\end{definition}

\begin{theorem}
A splitting surface $S$ for $B$ determines a unique lamination $\mathcal{L}_S$
fully carried by $B$. Conversely, any lamination $\mathcal{L}$ fully
carried by $B$ determines a splitting surface $S$, and $\mathcal{L}_S$
is a canonical blow-down of $\mathcal{L}$.
\end{theorem}
\begin{proof}
For each $I$-fiber $J$ of $\mathcal{V}$,   $J-\mathcal{N}(S)$ is a closed
subset. If $J-\mathcal{N}(S)$ has non-empty interior, then consider
the quotient $J'$ where each open interval of $int(J-\mathcal{N}(S))$ is
identified to a point. Doing this for each $I$-fiber, we get a quotient
complex $N(B)'$ with fibration $\mathcal{V}'$ over the branched surface
$B$, where the fibration is by points over the subcomplex of faces
of $B$ which do not carry $S$, and is a fibration by intervals over
$\pi(S)$. The complement of $\mathcal{N}(S)$ in $N(B)'$ is a
lamination $\mathcal{L}_S$ fully carried by $B$. Moreover,
$\mathcal{L}_S$ is unique
up to normal isotopy.

Conversely, suppose we have a lamination $\mathcal{L}$ fully
carried by $B$. Isotope $\mathcal{L}$ so that it meets the
endpoints of each fiber of $\mathcal{V}$, possibly collapsing
intervals of $\mathcal{V}$ to points over the faces of $B$ which
carry $\mathcal{L}$ with multiplicity one. Then the complementary
regions which meet $\partial_v N(B)$ consist of $I$-bundles over
an embedded surface $S$, which we may choose to intersect
$\partial_v N(B)$ in $p(B)$, so that the region is
$\mathcal{N}(S)$ in the topology induced from $M_\mathcal{L}$
using the path metric. Thus $S$ is a splitting surface. Consider
the intersection with each $I$-fiber $J$, and as before quotient
the open intervals of $J-\mathcal{N}(S)$ to points. This quotient
map gives a blow-down of $\mathcal{L}$ to $\mathcal{L}_S$.
\end{proof}

The branched surface $B$ has a cell structure $\tau$. A splitting
surface $S$ has an induced cell structure $\tau_S$ by pulling back
the cell structure $\tau$ under the quotient map $\pi$ (which
collapses each fiber of $\mathcal{V}$ to a point).
\begin{definition}
 A {\it splitting complex} is an embedded
complex $c$ in $N(B)$ whose cells are transverse to $\mathcal{V}$
and which has a cell structure $\tau_c$ induced by pullback by
$\pi$ from $\tau$. Moreover, $c$ is locally embedded in a surface
transverse to $\mathcal{V}$, in particular $int(c)$ is a surface.
Let $\partial c$ be the graph of all 1-cells of $\tau_c$ which are
not incident with two 2-cells of $\tau_c$, ie $\partial c=
c-int(c)$ . Then we also require that $\partial c\cap \partial_v
N(B) = p(B)$.
\end{definition}

A $B$-isotopy of a splitting complex $c$ is an isotopy of $c$
through splitting complexes which preserves intersections with
fibers of $\mathcal{V}$.
 So a splitting surface $S$
is also a splitting complex, and one should think of a splitting
complex as a subcomplex of a splitting surface. If $b\subset c$ is
a subcomplex of a splitting complex $c$, then let $B_1(b,c)$
consist of the subcomplex of $c$ by adjoining all closed cells
incident with $b$ to $b$. This is akin to taking the combinatorial
ball of radius 1 about $b$. Define inductively
$B_k(b,c)=B_1(B_{k-1}(b,c),c)$. Let the {\it radius} of a
splitting complex $c$ be $R(c)=\max \{k|B_k(p(B),c)\cap \partial
c=p(B)\}$. The radius measures the minimal combinatorial distance
from $p(B)$ to $\partial c-p(B)$.  Suppose that each cell of $B$
has at most $r$ cells (counted with multiplicity) incident with
it. Then for any splitting complex $c$, $B_k(p(B),c)$ has at most
$r^k |p(B)^{(0)}|$ cells.

\begin{lemma} \label{split}
There are only finitely many splitting complexes of $B$ with $\leq
N$ cells, up to $B$-isotopy, and there is an algorithm to
enumerate them.
\end{lemma}
\begin{proof}
The $B$-isotopy class of each splitting complex $c$ is determined
by the number of cells of each type and how they are glued
together along $\pi^{-1}(\tau^{(1)})$. To describe a splitting
complex with $\leq N$ 2-cells up to $B$-isotopy, we
 list for each 2-cell $f$ of $\tau$, the number of 2-cells
 of $c$ which map to $f$ under $\pi$. These 2-cells are numbered
 by the vertical order they come in the $I$-fibration of $\pi^{-1}(f)$.
 Then we give a pairing of a subset of edges of these cells, and a
 pairing of some edges with edges of $p(B)$.
 There are
 checkable restrictions on the pairings so that the quotient
 complex is a splitting complex. The pairings must preserve
 the orderings of the cells, and must not induce branching
 about $\pi^{-1}(\tau^{(0)})$, and there must be exactly
 one edge paired with each edge of $p(B)$.
 \end{proof}
 \begin{theorem}
 There is a sequence of splitting complexes $c_k$ of $B$ such that the radius
 $R(c_k)\to \infty$ as $k\to \infty$ if and only if there is
 a splitting surface $S$ of $B$.
 \end{theorem}
 \begin{proof}
 If there is a splitting surface $S$, then $R(B_k(p(B),S))=k-1 \to \infty$
 as $k\to \infty$, so we may take $c_k = B_k(p(B),S)$.

 Conversely, given a sequence of splitting complexes $c_k$, with
 $R(c_k)\to \infty$, we may assume that $\{R(c_k)\}$ is an
 increasing sequence. Then for $k\geq K$, the sequence of
 complexes $\{B_{R(c_K)}(p(B),c_k)\}$ must have a constant
 subsequence, since there are at most $|p(B)^{(0)}|r^{R(c_K)}$
 cells in these complexes, and thus there are only finitely
 many complexes up to $B$-isotopy, by lemma \ref{split}.
 Thus, we may find a subsequence $\{c_{k_i}\}$ such that for
 all $j>i$, $B_{R(c_{k_i})}(p(B),c_{k_i})$ is $B$-isotopic to $B_{R(c_{k_i})}(p(B),c_{k_j})$.
 We may $B$-isotope $B_{R(c_{k_j})}(p(B),c_{k_j})$ so that
 $B_{R(c_{k_i})}(p(B),c_{k_i}) =$$ B_{R(c_{k_i})}(p(B),c_{k_j})$.
 Moreover, we may choose a small $I$-bundle neighborhood of
 each $\{B_{R(c_{k_i})}(p(B),c_{k_i})\}$ so that the later stages remain disjoint from it.
 In the limit, the union of $B_{R(c_{k_i})}(p(B),c_{k_i})$ forms
 a splitting surface $S$.
\end{proof}

\begin{theorem} \label{lamalg2}
Given a branched surface $B$, there is a procedure which
terminates if and only if $B$ does not carry a lamination.
\end{theorem}
\rk{Remark}The referee points out that this theorem was already known
to Oertel.
\begin{proof}
 The algorithm proceeds by enumerating all splitting
complexes $c$ using $\leq |p(B)^{(0)}| r^N$ cells, by the
algorithm described in lemma \ref{split}, for $N=1,2,3,\ldots$. For
each $N$, one finds whether there is such
 a $c$ with $R(c)\geq N-1$, which is algorithmic since
 given $c$ there is a simple algorithm to compute $R(c)$.
 If no such complex $c$ exists,
 then $B$ carries no lamination fully, since
 if there were a splitting surface $S$, then $B_N(p(B),S)$
 would be a splitting complex with at most $\leq |p(B)^{(0)}| r^N$ cells
 and with radius $\geq N-1$.
\end{proof}

\section{Constructing essential branched surfaces}
Suppose we split a certain branched surface $B$ along a splitting
surface described before.  The inverse limit of an infinite
splitting is a lamination carried by $B$ by Lemma 4.2 of
\cite{GO}. Thus, the process of such splitting will stop unless
$B$ fully carries a lamination.  By Theorem 1 of \cite{GO}, if an
incompressible Reebless branched surface fully carries a
lamination, then it is an essential lamination. Therefore, in
order to make the algorithm work, we need to find finitely many
incompressible Reebless branched surfaces such that one of them
fully carries a lamination if $M$ is laminar.

We assume that $M$ is irreducible and has a triangulation $\tau$.
We may also assume that $M$ is atoroidal, otherwise $M$ contains
an essential torus, and thus has an essential lamination. By
\cite{JR}, there is an algorithm to find an efficient
triangulation $\tau$ of $M$, which has the property that the only
normal 2-sphere is the vertex linking 2-sphere, and there are only
finitely many normal tori (each bounding a solid torus).  It has
been shown in \cite{Br, Ga8} that if $M$ contains an essential
lamination, then $M$ contains a normal essential lamination with
respect to any given triangulation.  By putting normal disks
together, as on p. 122 of \cite{FO}, we can construct finitely
many branched surfaces such that any normal lamination is fully
carried by one of them. In fact, for any normal essential
lamination, by identifying all the normal disks of the same type
to one normal disk, we get a branched surface $B$ that is among
the finitely many branched surfaces constructed above and $B$
fully carries the essential lamination. By solving certain systems
of branch equations, we can algorithmically find all minimal
weight disks of contact in $B$ and split $B$ to eliminate all of
these disks of contact. This can be done in finitely many steps.
The argument of Claims 1--3, pages 122--123 of \cite{FO} generalizes
from incompressible surfaces to laminations to show that, once we
have split along all minimal weight disks of contact, we will have
a branched surface carrying a lamination isotopic to $\lambda$
which does not carry any 2-sphere, since $\tau$ was assumed to be
an efficient triangulation.

We say that a surface $S$ is {\it carried} by $N(B)$ if $S$ is
transverse to the $I$-fibers of $N(B)$.  In this paper, we always
assume our surface $S$ above is embedded in $N(B)$.  We define the
{\it weight} of a surface $S$  to be the number of intersection
points of $S$ with the 1-skeleton of the triangulation, and define
the {\it length} of an arc to be the number of intersection points
of this arc with the 2-skeleton.  Let $p_0$ and $p_1$ be surfaces
or arcs in $N(B)$ transverse to the $I$-fibers of $N(B)$.  We say
$p_0$ and $p_1$ are $B$-{\it parallel} if there is an isotopy $H:
F\times [0,1]\to N(B)$ such that $H(F\times i)=p_i$ ($i=0,1$) and
$H(q\times [0,1])$ is a subarc of an $I$-fiber of $N(B)$ for any
$q\in F$. We also call this isotopy a $B$-{\it isotopy}.

\begin{lemma}\label{L:torus}
Suppose $B$ fully carries a nowhere dense essential lamination
$\lambda$ and $B$ contains no disk of contact and carries no
2-spheres. Let $T$ be a torus carried by $N(B)$ and bounding a
solid torus $V$ in $M$. Then, after $K$ steps of splitting, where
$K$ depends on $B$ and $T$ (not on $\lambda$), $B$ can be split
into a branched surface that still carries $\lambda$ (up to
isotopy) but does not carry $T$.
\end{lemma}
\begin{proof}
If there are subarcs of $I$-fibers of $N(B)$ properly embedded in
$M-int(V)$, then the union of these subarcs is an $I$-bundle over
a compact surface whose horizontal boundary lies in $T$.  Let
$H\subset T$ be the horizontal boundary of this $I$-bundle.  Note
that we can assume $H$ consists of planar surfaces, otherwise,
since $M$ is irreducible, $M$ must be the union of $V$ and a
twisted $I$-bundle over a Klein bottle, which contradicts that $M$
contains an essential lamination (because such 3-manifolds are
double covered by lens spaces).  Then, by splitting along a union
of surfaces (in this $I$-bundle) with total weight bounded by the
weight of $H$, we get a branched surface that still carries $T$
and $\lambda$ but there is no subarc of any $I$-fiber properly
embedded in $M-int(V)$.

So we may assume that there is no subarc of any $I$-fiber of
$N(B)$ properly embedded in $M-int(V)$.  In the discussion below,
we always assume $T$ intersects $\lambda$ transversely.  We can
perform $B$-isotopy on parts of $\lambda$ pushing $\lambda$ into
$V$ as much as one can in the sense that after the $B$-isotopy,
for any point $x\in T\cap\lambda$, there is an arc
$\alpha\subset\lambda-int(V)$ such that, $x\in\alpha$, $\alpha$ is
not $B$-parallel to any arc in $T$, and the combinatorial length
of $\alpha$ is bounded by $K_1$ , where $K_1$ depends on $B$ and
$T$ (not on $\lambda$).

If $\lambda\cap T$ contains circles homotopically trivial in $T$,
then as in \cite{Br3}, there are trivial circles bounding disks in
$V\cap\lambda$ or $\lambda-int(V)$.  Then, after some standard
cutting and pasting on $\lambda$, we can eliminate the trivial
circles in $T\cap\lambda$.  Note that one can perform such cutting
and pasting on laminations because of the Reeb stability theorem
(see the proof of Lemma 2.1 of \cite{Br3}).  Moreover, since
$N(B)$ does not carry any 2-sphere, after such cutting and
pasting, the lamination is still transverse to the $I$-fibers of
$N(B)$.  After this operation, by Lemma 2.1 and Theorem 3.1 in
\cite{Br3}, $V\cap\lambda$ is either a union of meridional disks
or contains a sublamination by annuli (with at most one M\"{o}bius
band). Moreover, if the second case happens, any non annulus (or
M\"{o}bius band) leaf in $V\cap\lambda$ is non-compact and simply
connected, and each annulus leaf is $\pi_1$-injective in $V$.

Suppose $V\cap\lambda$ is a union of meridional disks.  In
general, a meridional disk (with fixed boundary) can wrap around
$T$ many times similar to a Reeb foliation.  However, we can
perform an isotopy (in fact a Dehn twist) near $T$ to unwrap these
disks, see Figure 6.1 in \cite{Ga8} for a schematic picture. Thus,
after this unwrapping, there must exist a point $x$ in
$T\cap\lambda$ and an arc $\beta\subset V\cap\lambda$ such that
$x\in\beta$, $\beta$ is not $B$-parallel to any arc in $T$, and
the length of $\beta$ is bounded by $K_2$, where $K_2$ depends on
$V\cap N(B)$ (not on $\lambda$).  By connecting $\beta$ and
$\alpha$ above together, we get an arc (puncturing through $T$)
with length less than $K_1+K_2$.  Hence, by splitting $B$ along
some disk with diameter less than $K_1+K_2$, one gets a branched
surface that carries a lamination isotopic to $\lambda$ but does
not carry $T$.

Now, we suppose $V\cap\lambda$ has a sublamination by
$\pi_1$-injective annuli (with at most one M\"{o}bius band). The
annular leaves that are $B$-parallel to subannuli in $T$ (fixing
the boundary) form a sublamination of $V\cap\lambda$. After some
$B$-isotopy, we can push these annuli (and the simply connected
leaves in between) out of $V$ (this can also be done by performing
some cutting and pasting on annuli in $T$).  Such a $B$-isotopy
may create some new intersection points in $T\cap\lambda$ and
these new intersection points lie in the simply connected leaves
of $V\cap\lambda$. After perturbing the simply connected leaves in
$V\cap\lambda$, for each new intersection point $x$ in
$T\cap\lambda$ (created during the $B$-isotopy above), there is an
arc $\alpha_x\subset\lambda$ connecting $x$ to an intersection
point that is fixed during the $B$-isotopy above. Moreover,
$\alpha_x$ is $B$-parallel to an arc in an annulus leaf that is
pushed out of $V$, and the length of $\alpha_x$ is bounded by a
number depending only on $T$.  So, for each point $x$ in
$T\cap\lambda$ after this $B$-isotopy, there is still an arc
$\alpha$ with length bounded by a number $K_1$ as above,
containing $x$, and not $B$-parallel to any arc in $T$.

So, we can assume $T\cap\lambda$ contains no annulus leaf that is
$B$-parallel to a subannulus in $T$, and $T\cap\lambda$ is not a
union of meridional disks.  In general, an annulus leaf in
$T\cap\lambda$ (with fixed boundary) can form a $monogon\times
S^1$ and wrap around $T$ many times.  However, after a Dehn twist
near $T$, any annulus is isotopic (fixing the boundary) to one
with smaller weight, as in the case of meridional disks above. In
particular, after this isotopy, there exists a point $x$ in the
boundary of this annulus leaf $A$ and an arc $\beta\subset A$ with
$x\in\partial\beta$ such that $\beta$ is not $B$-parallel to any
arc in $T$ and the length of $\beta$ is bounded by $K_2$, where
$K_2$ is as above.  Then, as before, by connecting $\alpha$ and
$\beta$, we get an arc (puncturing through $T$) with length less
than $K_1+K_2$.  Hence, by splitting $B$ along some disk with
diameter less than $K_1+K_2$, one gets a branched surface carrying
$\lambda$ but not $T$.

Therefore, in any of the cases, after splitting along a union of
surfaces with bounded total weight, we can get a branched surface
that carries $\lambda$ (up to isotopy) but does not carry $T$.
Since the total weight is bounded by a number that does not depend
on $\lambda$, we can enumerate all possible surfaces along which
we perform splitting as above, and after splitting we get finitely
many branched surfaces from $B$, one of which carries $\lambda$
but does not carry $T$.
\end{proof}
Note that after the splitting performed in the proof of Lemma
\ref{L:torus}, new disks of contact may appear, but we can always
find and eliminate them by another splitting and taking
sub-branched surfaces.

\begin{proposition}\label{P:monogon}
Suppose $B$ fully carries an essential lamination and $B$ contains
no disk of contact. Let $C$ be a component of $M-int(N(B))$ and
suppose $C$ contains a monogon.  Then, $C$ must be a solid torus
in the form of $D\times S^1$, where $D$ is a monogon.
\end{proposition}
\begin{proof}
Since $B$ fully carries an essential lamination and $B$ contains
no disk of contact, $\partial_hN(B)$ is incompressible in $M$. Let
$D$ be a monogon in $C$, ie, the disk $D$ is properly
embedded in $C$,  $\partial D$ consists of two arcs
$\alpha\subset\partial_vN(B)$ and $\beta\subset\partial_hN(B)$,
and $\alpha$ is a vertical arc in $\partial_vN(B)$.  Let $v$ be
the component of $\partial_vN(B)$ containing $\alpha$ and $N(v\cup
D)$ be a small regular neighborhood of $v\cup D$ in $M$.  Then,
the intersection of $N(v\cup D)$ and the component of
$\partial_hN(B)$ containing $\beta$ is a circle $\gamma$.  By the
construction, $\gamma$ is a trivial curve in $N(v\cup D)$.  Since
$\partial_hN(B)$ is incompressible, $\gamma$ must bound a disk in
$\partial_hN(B)$.   Thus, the component of $\partial_hN(B)$
containing $\beta$ is an annulus, $\partial C$ is a compressible
torus, and $D$ is a compressing disk for $\partial C$. Since $C$
is irreducible, $C$ must be a solid torus in the form of $D\times
S^1$, where $D$ is the monogon above.
\end{proof}

Suppose $A$ is an annulus carried by $N(B)$.  We say $A$ is a {\it
splitting annulus} if $\partial A$ lies in distinct components of
$\partial_vN(B)$.

\begin{proposition}\label{P:annulus}
Let $B$ be a branched surface fully carrying an essential lamination $\lambda$, and suppose $B$ has a monogon.  Then, there is a splitting annulus $A$ in $N(B)$ such that $A\cap\lambda=\emptyset$.
\end{proposition}
\begin{proof}
We can suppose $B$ does not have any disk of contact.  Let $C$ be
a component of $M-int(N(B))$ containing a monogon $D$.  Suppose
$\partial D=\alpha\cup\beta$, where $\alpha$ is a vertical arc in
$\partial_vN(B)$ and $\beta\subset\partial_hN(B)$.  By Proposition
\ref{P:monogon}, $C=D\times S^1$.  Let $\nu=\alpha\times S^1$ be
the corresponding component of $\partial_vN(B)$.  By the
end-incompressibility of the essential lamination,  we can split
$N(B)$ along $\lambda$ (by drilling a tunnel) and connect $C$ to
another component $W$ of $M-int(N(B))$, ie, there is a
vertical rectangle $R$ in $N(B)$ connecting $\nu$ and another
component of $\partial_v N(B)$ in $\partial W$ such that
$R\cap\lambda=\emptyset$. Note that if $R$ connects $\nu$ to
itself, one gets a compressing disk that is the union of $R$ and
two monogons, which contradicts that $\lambda$ is essential unless
$R$ is parallel to $\nu$.  So, existence of such a vertical
rectangle $R$ is guaranteed by the end-incompressibility of
$\lambda$. After removing a small neighborhood of $R$,  $N(B)$
becomes $N(B')$ which is a fibered neighborhood of another
branched surface $B'$ carrying  $\lambda$, and now $C$ and $W$
(connected through $R$) become a component $C'$ of $M-int(N(B'))$,
and $\nu$ becomes $\nu'$ which is a component of $\partial_v
N(B')$. Moreover, since $\lambda$ is end-incompressible, we can
assume that $W$ is not a $D^2\times I$ region.

Since $C=D\times S^1$, after this splitting above, $C'$ still
contains a monogon.  As $W$ is not a $D^2\times I$ region, $C'$
cannot be in the form of $monogon\times S^1$. Therefore, by
Proposition~\ref{P:monogon}, $B'$ cannot be incompressible. Since
$B'$ fully carries $\lambda$ and $B$ contains no disks of contact,
$B'$ must contain a disk of contact whose boundary must lie in
$\nu'$, ie, there must be a disk in $N(B')$ transverse to
the $I$-fibers and with boundary in $\nu'$.  By our construction
of $B'$ and $\nu'$, there must be a splitting annulus $A$ in
$N(B)$ connecting $\nu$ and another component of $\partial_vN(B)$
with $A\cap\lambda=\emptyset$ and $A\cap R$ a nontrivial arc in
$A$.
\end{proof}

Intuitively, if $M$ is atoroidal, one should be able to eliminate
all tori carried by a branched surface in finitely many steps of
splitting as in Lemma~\ref{L:torus}.  However,  the situation can
be very complicated if the branched surface carries infinitely
many tori. These tori can tangle together in a complicated way,
and it is not clear to us whether there is a simple way to deal
with it.  This problem can be simplified to a great extent, if we
use a special kind of triangulation, namely the one-efficient
triangulation due to Jaco and Rubinstein \cite{JR}. A
triangulation is one-efficient if every normal torus is either
thin or thick edge linking.  A trivial consequence of using a
one-efficient triangulation is that there are only finitely many
normal tori, which is the only thing we need for our purpose.

\begin{lemma}\label{L:splitting}
Suppose $M$ is an atoroidal 3-manifold with a one-efficient
triangulation $\tau$. Then there is an algorithm to construct
finitely many incompressible Reebless branch\-ed surfaces such
that every (nowhere dense) essential lamination normal with
respect to $\tau$ is fully carried by one of these branched
surfaces.
\end{lemma}
\begin{proof}

As before, we start with finitely many normal branched surfaces.
Since $\tau$ is one-efficient, any normal 2-sphere is normal
isotopic to the vertex linking 2-sphere. If a normal branched
surface $B$ carried a 2-sphere, then some component of $\partial_h
N(B)$ would be a 2-sphere parallel to the vertex linking 2-sphere,
so $B$ could not fully carry an essential lamination, since then
it would have a sphere leaf. So we may assume $B$ carries no
2-spheres. We eliminate all disks of contact, then take
sub-branched surfaces. Since there are only finitely many normal
tori, by Lemma~\ref{L:torus}, we can split these branched surfaces
in finitely many steps to construct finitely many branched
surfaces that do not carry any tori, and each normal essential
lamination is carried by one of them, then we get rid of disks of
contact again. So, these branched surfaces contain no Reeb
components. By taking sub branched surfaces if necessary, we can
assume each normal essential lamination is fully carried by one of
finitely many such branched surfaces.

Let $B$ be a branched surface constructed above.  $B$ contains no
disk of contact and does not carry any torus.  Suppose $B$ fully
carries an essential lamination $\lambda$ but $M-B$ has a monogon.
By Proposition~\ref{P:annulus}, there is a splitting annulus $A$
($A\cap\lambda=\emptyset$) connecting two components $v_1$ and
$v_2$ of $\partial_vN(B)$.  Next, we analyze annuli carried by
$N(B)$ with exactly two boundary circles lying in $v_1$ and $v_2$
respectively by solving a system of branch equations as in
\cite{FO}.  Note that we are considering surfaces with boundary,
so some equations are like $x_i+x_j=x_k-1$ and this system of
linear equations is not homogeneous. As in \cite{FO}, there is a
one-to-one correspondence between non-negative integer solutions
to this system of branch equations and surfaces carried by $B$
with exact two boundary circles lying in $v_1$ and $v_2$
respectively.  Such a solution is a point in $\mathbb{R}^m$ with
non-negative integer coordinates, where $m$ is the number of
variables.

Suppose $A_1=(a_1,\dots,a_m)$ and $A_2=(b_1,\dots,b_m)$ are two
solutions to the system above, and suppose
$\chi(A_1)=\chi(A_2)=0$, where $\chi(A_j)$ denotes the Euler
characteristic of the surface $A_j$.  If $a_i\le b_i$ for each
$i$, then $T=(b_1-a_1,\dots, b_m-a_m)$ is a non-negative integer
solution to the corresponding homogeneous system, and hence $T$ is
a union of closed surfaces carried by $B$.  Moreover, $A_1+T=A_2$
implies $\chi(A_1)+\chi(T)=\chi(A_2)$.  Hence, $\chi(T)=0$.  Since
$B$ does not carry any 2-sphere or torus, this is impossible.  So,
there are no such pairs of annuli among integer solutions.  Now,
we interpret zero Euler characteristic as a linear equation and
add this equation to the system above.  Then, every non-negative
integer solution to this new system gives us an annulus (or a pair
of M\"{o}bius bands).  Since there are no two annuli $A_1$ and
$A_2$ as above, there are only finitely many non-negative integer
solutions to this new system, and hence there are only finitely
many possible splitting annuli connecting $v_1$ and $v_2$.  Note
that since this is a system of linear equations with integer
coefficients, the solution space must be bounded in the region of
$x_i\ge 0$ (for all $i$), otherwise there would be infinitely many
non-negative integer solutions.  Thus, one can calculate the
maximum for each coordinate, and list all non-negative integer
solutions, ie, one can list all possible splitting annuli.

After splitting along these splitting annuli, we can eliminate all
monogons and get finitely many branched surfaces.  By taking sub
branched surfaces if necessary, we have that, after isotopy, each
normal essential lamination is fully carried by one of these
branched surfaces.

Finally, for each of these branched surfaces, we  subdivide the
triangulation so that the branched surface lies in the 2-skeleton,
and then using the algorithm described in lemma \ref{incReeb} to
check that each branched surface is incompressible and Reebless.
\end{proof}
\begin{remark}
\begin{enumerate}
\item If $M$ is a small Seifert fiber space, one can always recognize this
manifold if it does not admit a one-vertex triangulation
\cite{JR}. Moreover, using layered solid tori, one can construct a
nice triangulation for a small Seifert fiber space \cite{J} that
also makes the proof work, though it may not be one-efficient.

\item One does not need a one-efficient triangulation to eliminate all
Reeb branched surfaces.  Suppose $T$ is a torus carried by $N(B)$
and bounding a solid torus $V$, and $V\cap N(B)$ fully carries a
sublamination of a Reeb foliation of a solid torus.  Then, each
$I$-fiber of $N(B)$ can intersect $T$ in at most 2 points.
Otherwise, since $T$ is separating, if an $I$-fiber intersects $T$
in more than 2 points, there must be a subarc of this $I$-fiber
properly embedded in the solid torus $V$. This contradicts a
well-known fact that there is no properly embedded compact arc in
a solid torus that is transverse to the Reeb foliation (extended
from the Reeb lamination).  Therefore, $B$ contains only finitely
many Reeb components.

\item By replacing a leaf of a lamination by an $I$-bundle over this
leaf and then deleting the interior of this $I$-bundle, one can
change every lamination to a nowhere dense one (see Remark 4.4,
\cite{GO}). Gabai proved that, up to isotopy, every nowhere dense
essential lamination in an atoroidal 3-manifold is fully carried
by one of finitely many essential branched surfaces (Theorem 6.5
of \cite{Ga8}). However, step 1 of the proof of Theorem 6.5 of
\cite{Ga8} uses Plante's argument on the limit of normal annuli,
and it does not clearly give an algorithm to find these finitely
many essential branched surfaces.  In particular, if a branched
surface carries infinitely many tori, the picture of those normal
annuli can be very complicated, and it is not clear how to
algorithmically analyze the limit of these annuli.
\end{enumerate}
\end{remark}

\begin{theorem} \label{alg}
There is an algorithm to decide whether a 3-manifold contains an essential lamination.
\end{theorem}
\begin{proof}
The first step of the algorithm is to modify the triangulation of
$M$ to a one-efficient triangulation $\tau$ using the algorithm in
\cite{JR}.  By \cite{JR}, in finitely many steps, we either get a
one-efficient triangulation, or $M$ contains an incompressible
torus, or we can recognize $M$ as certain small Seifert fiber
space.  An incompressible torus can be found, if one exists, by
\cite{JO} or by Algorithm 8.2 of \cite{JT}, and essential
laminations in small Seifert fiber spaces are classified by
\cite{Br3, EHN, JN, Na}. So, we know whether $M$ is laminar in
these exceptional cases, and we can assume our triangulation for
$M$ is one-efficient.

By Lemma~\ref{L:splitting}, we can algorithmically construct
finitely many incompressible Reebless branched surfaces such that
every nowhere dense essential lamination normal with respect to
$\tau$ is fully carried by one of them.  Then, if one of these
incompressible Reebless branched surfaces fully carries a
lamination (hence it is an essential lamination), the procedure
described in theorem \ref{lamalg1} will stop, which means $M$ is
laminar. If none of them fully carries a lamination, the procedure
described in theorem \ref{lamalg2} will stop, which means $M$
contains no essential lamination.
\end{proof}

\section{Recognizing Reebless foliations}
In this section, we will assume the reader is familiar with the
notion of a sutured manifold, introduced by Gabai in \cite{Ga2}.
We will follow the notation of \cite{Ga2}.

\rk{Remark} We note that if a manifold $M$ is atoroidal and
admits a Reebless foliation, then the foliation is taut. But in
the case of irreducible toroidal manifolds, there are examples
that admit Reebless foliations, but no taut foliation \cite{BNR}.
It seems to be unknown in general when a toroidal manifold admits
a taut foliation. Thus, in the theorems below, we will restrict
ourselves to considering Reebless foliations.

\begin{theorem} \label{taut}
 Given a sutured manifold
$(M,\g)$ with triangulation $\tau$ and $\g\subset \tau^{(1)}$,
 there is an algorithm to determine if $(M,\g)$ is
taut.
\end{theorem}
\begin{proof}
By theorem 4.2 of \cite{Ga2}, a connected taut sutured manifold
$(M, \gamma)$ has a sutured manifold hierarchy
$$(M_0,\g_0) \overset{S_1}{\longrightarrow} (M_1,\g_1) \overset{S_2}{\longrightarrow}
\cdots \overset{S_n}{\longrightarrow} (M_n,\g_n)$$ such that
$(M_n,A(\g_n))= (R\times I,\partial R\times I)$ and
$R_+(\g_n)=R\times 1$. Moreover, for each component $V$ of
$R(\g_i)$, $S_{i+1}\cap V$ is a union of $k (\geq 0)$ parallel
oriented non-separating simple closed curves or arcs. Moreover, if
$(M,\g)$ has a sutured manifold hierarchy such that no component
of $R(\g_i)$ is a compressible torus, then $(M,\g)$ is taut by
corollary 5.3 of \cite{Ga2}.

The algorithm to determine if $(M,\g)$ is taut proceeds by running
two procedures. The first procedure constructs subdivisions $\nu$
of the triangulation $\tau$ searching for surfaces $(S,\partial S)
\subset (M,\g) \cap \nu^{(2)}$ such that $[S,\partial S] =
[R_+,\g] \in H_2(M,\g)$, and $x(S)<x(R_{\pm})$. If $(M,\g)$ is not
taut, then this procedure terminates in a finite number of steps.

The second algorithm searches for a sutured manifold hierarchy.
Again, it constructs subdivisions $\nu$ of $\tau$, and searches
for oriented surfaces $S_i\subset \nu^{(2)}$ such that they
satisfy the conditions for a sutured manifold hierarchy. By
Algorithm 9.7 of \cite{JT}, there is an algorithm to check that
$(M_n,A(\g_n))\cong (R\times I,
\partial R\times I)$. There are also algorithms to check the
other conditions listed above. By Gabai's theorem, if $(M,\g)$ is
taut, then this procedure will halt in finitely many steps by
finding a sutured manifold hierarchy which certifies that $(M,\g)$
is taut.
\end{proof}

\begin{theorem}
There is an algorithm to decide whether an orientable 3-man\-if\-old
contains a Reebless foliation.
\end{theorem}
\begin{proof}
First, we describe the algorithm.

{\bf Step 1}\qua We check that $M$ is irreducible using an algorithm of
Rubinstein \cite{Ru2} and Thompson \cite{Tho}. If it is not
irreducible, then we check to see if $M$ is $S^2\times S^1$, by
cutting along a non-separating sphere, cap off the resulting
boundary with balls, and check if we have a ball again.  Otherwise
$M$ is not prime, and $M$ cannot have a Reebless foliation. If $M$
is irreducible then we continue.

{\bf Step 2}\qua We check to see if $M$ is toroidal. If $M$ is toroidal and
irreducible, then it has a transversely orientable Reebless
foliation by the method of Corollary 6.5 of \cite{Ga2}.

{\bf Step 3}\qua If $M$ is irreducible and atoroidal, then we construct
finitely many normal essential branched surfaces $B$ which fully
carry every normal essential lamination using Lemma
\ref{L:splitting}. We check to see if $B$ carries a lamination. We
triangulate the complementary regions of $B$. Then we use
algorithm 8.2 of \cite{JT} to find the maximal $I$-bundle
$(C,\partial_v C)\subset (M-int N(B),\partial_v N(B))$. We then
delete from $M-int N(B)$ unions of components $I$-bundle
components $C'\subset C$ in all $2^{\beta_0(C)}$ possible ways. We
use the algorithm described in theorem \ref{taut} to determine if
the regions of $M-int(N(B)\cup C')$ are taut sutured manifolds. If
all of them are, then by a theorem of Gabai \cite{Ga2}, $M-int
(N(B)\cup C')$ has foliations transverse to $\partial_v(N(B)\cup
C')$, with $\partial_hN(B)-C'$ as leaves. We may extend these
foliations and the lamination carried by the branched surface to a
Reebless foliation, using the method of the proof of theorem 5.1,
pages 471--477 in \cite{Ga2} (see also constructions 4.16 and 4.17,
pages 498--500 of \cite{Ga10}), since $B$ contains no disk of
contact.

To see that the algorithm works, we need to show that if $M$ has a
Reebless foliation, then for the branched surfaces constructed in
Lemma \ref{L:splitting}, one of the essential branched surfaces
union some $I$-bundles has as complementary regions taut sutured
manifolds which may be extended to a foliation. If $M$ has a
Reebless foliation, then by theorem 4.4 of \cite{Ga8}, there
exists a normal essential lamination such that each complementary
region of the normal essential branched surface $B'$ carrying this
lamination is either a taut sutured manifold or an $I$-bundle,
since the normal lamination is obtained from the foliation by
first splitting, then evacuating a taut sutured manifold.
Moreover, these regions are $\pi_1$-injective in the 3-manifold
group, by condition (ii) of definition 4.2 \cite{Ga8}. The normal
essential lamination fully carried by $B'$ is also fully carried
by one of the branched surfaces $B$ constructed in Lemma
\ref{L:splitting}, and thus there is an essential normal branched
surface $B_0$ that is a splitting of both $B$ and $B'$.  So, there
is a union of product regions $\mathcal{C}$ such that
$N(B_0)\cup\mathcal{C}=N(B)$. Let $\lambda$ be the normal
essential lamination fully carried by $B_0$, $B$ and $B'$, and
suppose $\lambda$ extends to a Reebless foliation $\mathcal{F}$.
For any $I$-bundle $E\in\mathcal{C}$, we suppose that the
horizontal boundary of $E$ lies in $\lambda$ and the vertical
boundary of $E$ consists of essential annuli in $M-N(B_0)$.
Moreover, since $B_0$ is a splitting of $B$ and $B$ contains no
splitting annulus as in the construction in Lemma
\ref{L:splitting}, we also assume that $E$ is not an
$annulus\times I$ or an $I$-bundle over a M\"{o}bius band. Next,
we show that after isotopies, $\mathcal{F}$ is transverse to the
$I$-fibers of $E$ for any $E$.  Note that the case that some
vertical boundary component of $E$ is not an essential annulus
(ie, bounds a $D^2\times I$ region) and the case that $E$ is an
$I$-bundle over an annulus or a M\"{o}bius band with a vertical
boundary component in $\partial_vN(B_0)$ are easy to prove by
similar arguments.

We can split $\mathcal{F}$ into a lamination $\mathcal{F}'$ so
that $\lambda$ is a sublamination of $\mathcal{F}'$ and,  after
isotopy, $\mathcal{F}'\cap E$ is an essential lamination  in $E$.
By our assumption that $E$ is not an $I$-bundle over   an annulus
or a M\"{o}bius band, the double of $E$ is an   $I$-bundle over a
closed surface of genus greater than 1.    A theorem of Brittenham
\cite{Br97} says that any essential    lamination in an $I$-bundle
over a closed surface with    genus $>1$, containing the boundary
as leaves, can be    isotoped to be transverse to the $I$-fibers.
Thus, we can assume $\mathcal{F}'\cap E$ is transverse to the
$I$-fibers    of $E$.  Let $N(B_0)\cup E=N(B_1)$ be a fibered
neighborhood     of another essential branched surface $B_1$.
Since the     vertical boundary of $E$ consists of essential
annuli and     $E$ is not an $I$-bundle over an annulus or a
M\"{o}bius band,      $B_1$ contains no splitting annulus that
cuts through $E$.      Note that by the construction in
\cite{Ga8}, every $I$-fiber  of $N(B')$ (and hence $N(B_0)$) is
transverse to the       foliation $\mathcal{F}$.  If there is an
$I$-fiber of       $E$ that cannot be transverse to the foliation
$\mathcal{F}$, then there must be an arc $\alpha$ in a leaf of
$\mathcal{F}'$ and a subarc $\beta$ of an $I$-fiber in $\partial
E$ such that $\alpha\cup\beta$ bounds a monogon disk in
$M-\mathcal{F'}$. Similar to the proofs of Propositions
\ref{P:monogon} and \ref{P:annulus}, such a monogon disk implies
that $N(B_1)$ contains a splitting annulus cutting through $E$,
which gives a contradiction. Hence, after isotopies, the foliation
$\mathcal{F}$ is transvere to the $I$-fibers of $E$ for any
$E\in\mathcal{C}$.

Therefore, by our assumptions on $B'$, $\lambda$ and
$\mathcal{F}$, there must be a union of product regions $C'$ in
$M-int(N(B))$ such that $M-int(N(B)\cup C')$ consists of taut
sutured manifolds.  Conversely, if we find such a $C'$ for $B$, we
can conclude that the essential lamination fully carried by $B$
extends to a Reebless foliation by \cite{Ga2,Ga10} because the
complementary regions of $B$ are essential.  Thus, the algorithm
will succeed if and only if the manifold has a Reebless foliation.
\end{proof}

\section{Conclusion}
The algorithm we have described is unsatisfying, since the proof
that the algorithm terminates doesn't give us any idea of how long
the algorithm will take, and it seems that it would be nearly
impossible to implement on a computer. It would be interesting to
get an upper bound on how much one needs to split an essential
branched surface to get a laminar branched surface. If the bound
were good enough, then one might be able to answer the following
question.

\rk{Question}Is there an NP algorithm to determine if a 3-manifold is
laminar?\medskip

Other questions of interest:
\begin{enumerate}
\item
Is there a finite collection of laminar branched surfaces in $M$
which carry every essential lamination, and if so, is there an
algorithm to find them?
\item
Are there algorithms to determine if a manifold has a tight
lamination \cite{Br}, an $\RR$-covered foliation \cite{C}, a
slithering \cite{Th4}, or a pseudo-Anosov flow \cite{Fr} (these
references give definitions of these objects, not necessarily
original sources for the definitions)?
\item
Is there an algorithm to determine if the fundamental group of a
3-manifold has a circular ordering? If an atoroidal manifold $M$
has a taut transversely orientable foliation, then it has been
shown by Bill Thurston that $\pi_1(M)$ has a circular ordering,
also proven by Calegari and Dunfield in Theorem 6.2 of \cite{CD}.
\item
Is there an algorithm to find all essential laminations in Dehn
fillings on a link? What we have in mind here would be to describe
a collection of branched surfaces in a link complement, and
a description for each branched surface of which Dehn fillings
on the link have essential laminations which meet the link
complement in a lamination carried by that branched surface.
This has been done for the figure eight knot complement by
Tim Schwider \cite{Schw}.

\end{enumerate}

\end{document}

%% file: gtoutput.tex
\def\ifplaintex{\expandafter\ifx\csname documentclass\endcsname\relax}

\ifplaintex 
\else
\fi

\expandafter\ifx\csname beginpicture\endcsname\relax
\expandafter\ifx\csname documentclass\endcsname\relax
\input pictex \else
\input prepictex \input pictex \input postpictex \fi\fi

\def\gt{{\mathsurround=0pt\it $\cal G\mskip-2mu$eometry \&\ 
$\cal T\!\!$opology}}        

\def\gtp{{\mathsurround=0pt\it $\cal G\mskip-2mu$eometry \&\ 
$\cal T\!\!$opology $\cal P\!$ublications}}  

\def\lognumber#1{\def\thelognumber{#1}}
\def\volumenumber#1{\def\thevolumenumber{#1}}
\def\papernumber#1{\def\thepapernumber{#1}}
\def\volumeyear#1{\def\thevolumeyear{#1}}

\def\pagenumbers#1#2{\def\startpage{#1}\def\finishpage{#2}}
\def\published#1{\def\publishdate{#1}}
\def\proposed#1{\def\theproposer{#1}}
\def\seconded#1{\def\theseconders{#1}}
\def\received#1{\def\receiveddate{#1}}
\def\revised#1{\def\reviseddate{#1}}
\def\accepted#1{\def\accepteddate{#1}}

\def\asciiaddress#1{\def\theasciiaddress{#1}}
\def\asciiemail#1{\def\theasciiemail{#1}}

\let\\\par\let\thelognumber\relax
\let\thevolumenumber\relax\let\thepapernumber\relax
\let\thevolumeyear\relax\let\thesamplenumber\relax\let\startpage\relax
\let\finishpage\relax\let\publishdate\relax\let\receiveddate\relax
\let\reviseddate\relax\let\accepteddate\relax\let\theasciititle\relax
\let\theasciiauthors\relax\let\theasciiaddress\relax
\let\theasciiabstract\relax
\let\theasciiemail\relax\let\theshortauthors\relax\let\theshorttitle\relax

\long\def\maketitlep{   

\count0=\startpage

\gt\hfill      
\beginpicture
\setcoordinatesystem units <0.33truein, 0.33truein> point at 2.2 0.9
\setplotsymbol ({$\cal G$})
\plotsymbolspacing=9truept
\circulararc 315 degrees from 0 1 center at 0 0
\setplotsymbol ({$\cal T$})
\circulararc 315 degrees from 1 -1 center at 1 0
\endpicture
\break
{\small\ifx\thesamplenumber\relax 
Volume \else Sample
\fi\thevolumenumber\ (\thevolumeyear)
\startpage--\finishpage\nl
Published: \publishdate}
\vglue 0.5truein plus 0.4fil minus 0.1truein

{\parskip=0pt\leftskip 0pt plus 1fil\def\\{\par\smallskip}{\ifplaintex\large
\else\Large\fi\bf\thetitle}\par\medskip}   

\vglue 0pt plus 0.1fil 

{\parskip=0pt\leftskip 0pt plus 1fil\def\\{\par}{\sc\theauthors}
\par\medskip}

\vglue 0pt plus 0.1fil 

{\small\parskip=0pt\let\newline\\
{\leftskip 0pt plus 1fil\def\\{\par}{\sl\theaddress}\par}
\expandafter\ifx\theemail\relax    
\relax\else\vglue 5pt plus 0.02fil minus 2pt\def\\{\stdspace{\rm 
and}\stdspace} 
\cl{Email:\stdspace\tt\theemail}\fi
\ifx\theurl\relax                  
\relax\else\vglue 5pt plus 0.02fil minus 2pt\def\\{\stdspace{\rm 
and}\stdspace}
\cl{URL:\stdspace\tt\theurl}\fi\par}

\vglue 7pt plus 0.3fil minus 3pt

{\bf Abstract}
\vglue 5pt plus 0.1fil minus 2pt

\theabstract

\vglue 7pt plus 0.3fil minus 3pt

{\bf AMS Classification numbers}\quad Primary:\quad \theprimaryclass

Secondary:\quad \thesecondaryclass

\vglue 5pt plus 0.3fil minus 2pt

{\bf Keywords}\quad \thekeywords

\vglue 10pt plus 0.5fil minus 5pt

{\small  Proposed: \theproposer\hfill Received: \receiveddate\nl
Seconded: \theseconders\hfill 
\ifx\reviseddate\relax                         
Accepted: \accepteddate                        
\else
Revised: \reviseddate                          
\fi}
\eject
}       

\let\maketitlepage\maketitlep
\let\maketitle\maketitlepage

\font\phead=cmsl9 scaled 950
\font=cmsl9 scaled 1050
\font\pnum=cmbx10 scaled 913
\font\lnum=cmbx10 
\font\pfoot=cmsl9 scaled 950
\font=cmsl9 scaled 1050
\ifplaintex
\headline{\vbox to 0pt{\vskip -4.5mm\line{\small\phead\ifnum
\count0=\startpage ISSN 1364-0380 (on line)
1465-3060 (printed) \hfill {\pnum\folio}\else\ifodd\count0\def\\{ }%
\ifx\theshorttitle\relax\thetitle\else\theshorttitle\fi\hfill{\pnum\folio}
\else\def\\{ and }{\pnum\folio}\hfill\ifx\theshortauthors\relax\theauthors
\else\theshortauthors\fi\fi\fi}\vss}}
\footline{\vbox to 0pt{\vglue 0mm\line{\small\pfoot\ifnum\count0=\startpage
\copyright\ \gtp\hfill\else
\gt, Volume \thevolumenumber\ (\thevolumeyear)\hfill\fi}\vss
}}
\else
\makeatletter
\def\@oddhead{{\small\lhead\ifnum\count0=\startpage ISSN 1364-0380 (on line)
1465-3060 (printed) \hfill {\lnum\number\count0}\else\ifodd\count0
\def\\{ }\ifx\theshorttitle\relax \thetitle \else\theshorttitle\fi\hfill
{\lnum\number\count0}\else\def\\{ and }{\lnum\number\count0}
\hfill\ifx\theshortauthors\relax 
\theauthors\else\theshortauthors\fi\fi\fi}}\def\@evenhead{\@oddhead}
\def\@oddfoot{\small\lfoot\ifnum\count0=\startpage\copyright\ \gtp\hfill\else
\gt, Volume \thevolumenumber\ (\thevolumeyear)\hfill\fi}
\def\@evenfoot{\@oddfoot}
\makeatother
\fi

\newwrite\gtoutfile
\long\gdef\makeheadfile{  
{\def\\{, }\def\s{ }
\immediate\openout\gtoutfile head.xxx
\immediate\write\gtoutfile{To: math@arxiv.org}
\immediate\write\gtoutfile{Subject: put or rep NNNNN:pppp}
\immediate\write\gtoutfile{--text follows this line--}
\immediate\write\gtoutfile{Proxy-for: \ifx\theasciiauthors\relax
\theauthors\else\theasciiauthors\fi\s<\ifx\theasciiemail\relax\theemail\else\theasciiemail\fi>}
\immediate\write\gtoutfile{\noexpand\\}
\immediate\write\gtoutfile{Authors: \ifx\theasciiauthors\relax
\theauthors\else\theasciiauthors\fi}
{\def\\{ }\immediate\write\gtoutfile{Title: \ifx\theasciititle\relax
\thetitle\else\theasciititle\fi}}
\immediate\write\gtoutfile{Subj-class: GT or SG or MG etc}
\immediate\write\gtoutfile{MSC-class: \theprimaryclass\ifx\thesecondaryclass\relax\else, \thesecondaryclass\fi}
\immediate\write\gtoutfile{Journal-ref: Geom. Topol. \thevolumenumber
(\thevolumeyear) \startpage-\finishpage}
\immediate\write\gtoutfile{Comments: Published by Geometry and Topology at}
\immediate\write\gtoutfile{\s\s http://www.maths.warwick.ac.uk/gt/GTVol\thevolumenumber/paper\thepapernumber.abs.html}
\immediate\write\gtoutfile{\noexpand\\}
\immediate\write\gtoutfile{}
\ifx\theasciiabstract\relax
\immediate\write\gtoutfile{\theabstract}\else
\immediate\write\gtoutfile{\theasciiabstract}\fi
\immediate\write\gtoutfile{}
\immediate\write\gtoutfile{\noexpand\\}
\immediate\write\gtoutfile{}
\immediate\closeout\gtoutfile}}  

\def\maketitlepage{\maketitlep\makeheadfile}
\let\maketitle\maketitlepage

%% file: rlepsf.tex
\input epsf
\def\relabelbox{%
  \hbox\bgroup%
}%
\def\endrelabelbox{%
\egroup%
}%
\def\relabel #1#2 {%
  \special{ps:/a {} def}%
  \smash{\rlap{#2}}%
}%
\def\adjustrelabel <#1,#2> #3#4 {%
  \special{ps:/a {} def}%
  \smash{\rlap{\kern #1 \raise #2\hbox{#4}}}%
}%
\def\extralabel <#1,#2> #3 {\smash{\rlap{\kern #1 \raise #2\hbox{#3}}}}%